\documentclass[a4paper,10.5pt,twoside]{article}
\usepackage[top=3cm, bottom=3cm, inner=3cm, outer=3cm, includehead]{geometry}
\usepackage{fancyhdr}
\usepackage{xurl}
\usepackage{graphicx}
\usepackage{alltt}
\usepackage{amsmath}
\usepackage{amsthm}
\usepackage{amssymb}
\usepackage[hidelinks, pdftex]{hyperref}
\usepackage[T1]{fontenc}
\usepackage[utf8]{inputenc}
\usepackage{lmodern}
\usepackage{csquotes}
\usepackage[sorting=none]{biblatex}
\usepackage{xcolor}
\usepackage{multirow}
\usepackage{booktabs}
\usepackage[english]{babel}

\newtheorem{theorem}{Theorem}
\newtheorem{corollary}{Corollary}
\newtheorem{remark}{Remark}

\begin{document}
\fancyhead[LE]{\thepage\ \ \ \ LaRuez, Rooney}
\fancyhead[RO]{Draft, (2026)\ \ \ \ \thepage}

\begin{center}
\LARGE
\textbf{Preferential Attachment as a Simpliciality-Enforcing
Mechanism in Hypergraphs}\\[12pt]
\normalsize
\textbf{Jason LaRuez,\footnote{Rochester Institute of Technology,
jasonlaruez@gmail.com, corresponding author}
Brendan Rooney\footnote{Rochester Institute of Technology, brsma@rit.edu}}\\[4pt]
\end{center}

\begin{abstract}
\normalsize
Higher-order networks, represented as hypergraphs, enable direct modeling of multi-body interactions of arbitrary size. Hypergraph representations of real-world systems have been observed to exhibit high \emph{simpliciality} --- the tendency for subsets of hyperedges to also appear as hyperedges --- yet the generative mechanisms responsible for this structure are poorly understood. We introduce a generalized preferential attachment hypergraph model in which both hyperedge size $Y_t$ and the number of new nodes per step $X_t$ are drawn from arbitrary distributions, and derive analytically, using a mean-field approximate master equation approach, that the stationary hyperdegree distribution follows a power law whose exponent depends only on the ratio $p = E[X_t]/E[Y_t]$, independent of the shapes of the underlying distributions. Crucially, both $X_t$ and $Y_t$ can be estimated directly from any timestamped hypergraph dataset via a backward-stepping procedure, enabling the model to be fit without parametric assumptions. Applying a nonlinear extension of the model to eight real-world hypergraph datasets, we find that the simplicial fraction increases monotonically with the strength of preferential attachment up to the gelation transition at $\alpha > 1$, establishing preferential attachment as a simpliciality-enforcing mechanism.
\end{abstract}

\section{Introduction}\label{s:intro}

Networks have long served as modeling tools for systems with complex interaction patterns, including social dynamics \cite{OpinionFission}, epidemic spreading \cite{CoevolvingSISNetworkModel}, and technological infrastructure \cite{Infrastructure}. However, many real-world systems involve interactions among groups of arbitrary size that cannot be faithfully reduced to pairwise edges, motivating the use of \emph{hypergraphs} as modeling tools \cite{SimplicialClosure_LinkPrediction, XGI_data}.

A hypergraph $H = (V, E)$ consists of a node set $V$ and a collection $E$ of \emph{hyperedges}, each a subset of $V$ of size two or more. A \emph{simplicial complex} is a hypergraph with an additional structural property: for every hyperedge present. Given any hypergraph, one may construct its \emph{induced simplicial complex} by taking the downward closure of its edge set: for each hyperedge we add any of its subsets not present as a hyperedge. The downward closure property enables the application of tools from algebraic topology --- including Betti numbers and the Euler characteristic --- to the study of neural networks, social contagion, and synchronization \cite{PHNetworks, SimplicialContagion, SynchronizationHGSC, Brain}. However, hypergraph and simplicial complex models of the same system can exhibit qualitatively distinct dynamical and structural behavior \cite{SynchronizationHGSC, Encapsulation, Hyperedge_Overlap}, making it important to understand the conditions under which a hypergraph closely approximates its induced simplicial complex.

A key empirical observation is that real-world hypergraphs tend to be \emph{highly simplicial}: a large fraction of their hyperedges are already downwardly closed, so that the induced simplicial complex requires relatively few additions \cite{SimplicialClosure_LinkPrediction,Simpliciality, Simplicial_Pairs}. Despite significant interest in quantifying simpliciality \cite{Encapsulation,Simpliciality,Simplicial_Pairs}, the generative mechanisms producing simplicial structure in hypergraph data remain poorly understood.

The question of which generative rules produce scale-free degree distributions in graphs was answered by Barabási and Albert \cite{BA}, who showed that preferential attachment is a sufficient mechanism: new edges attach to existing nodes in proportion to their degree, producing ``rich-get-richer'' dynamics. Several preferential attachment models have since been proposed for hypergraphs \cite{EvolvingHGModel, PA_HG, PA_HG_LocalWorld, PA_HG_Poisson_Weighted}, but their relationship to simpliciality has not been studied. We address this gap with three contributions.

First, we introduce a generalized preferential attachment hypergraph model in which both the hyperedge size $Y_t$ and the number of new nodes per step $X_t$ are drawn from arbitrary distributions. Using a mean-field approximate master equation (AME) approach in the thermodynamic limit, we derive that the stationary hyperdegree distribution is a power law $P(d_H) \sim d_H^{-\gamma}$ with exponent $\gamma = 1/(1-p) + 1$ depending only on $p = E[X_t]/E[Y_t]$ --- regardless of the specific forms of the distributions of $X_t$ and $Y_t$ beyond their means. Second, we show that $X_t$ and $Y_t$ can be estimated directly from any timestamped hypergraph dataset via a simple backward-stepping procedure, grounding the model in empirical network statistics without parametric assumptions. Third, using a nonlinear extension of the model (NLPA), we demonstrate across eight real-world datasets that preferential attachment monotonically increases the simplicial fraction up to the gelation transition, establishing preferential attachment as a structurally coherent simpliciality-enforcing mechanism.

\section{Simpliciality}\label{s:simpliciality}

We quantify the similarity between a hypergraph and its induced simplicial complex using the \emph{simplicial fraction} $\sigma_{SF}$ introduced by Landry et al.~\cite{Simpliciality}. A hyperedge $e \in E$ is \emph{downwardly closed} if every non-empty proper subset of $e$ belongs to $E$. Letting $\hat{E} = \{e \in E \mid \mathcal{P}(e) \setminus \{\emptyset\} \subseteq E\}$ denote the set of downwardly closed hyperedges, the simplicial fraction is:
\begin{equation}
    \sigma_{SF} = \frac{|\hat{E}|}{|E|}.
\end{equation}
We have $\sigma_{SF} = 1$ if and only if $H$ is itself a simplicial complex. Following Landry et al. we restrict computations to hyperedges of size at least two. This is because many models do not allow for hyperedges containing singletons, thus total downward closure is too strict a requirement. Also, hyperedges of the minimum size are not counted towards the simplicial fraction since they are closed by default, and their inclusion can skew the simplicial fraction (thus hyperedges of size three or larger are checked for subsets of size two or larger).

The simplicial fraction can be computed efficiently via dynamic programming: a hyperedge $e$ is downwardly closed if and only if all of its $(|e|-1)$-sized subsets belong to $E$ and are themselves downwardly closed. Sorting hyperedges in ascending order of size and applying this observation inductively yields an $O(|E| \cdot \max_e |e|)$ algorithm that avoids the exponential memory cost of explicitly enumerating all missing subsets.

\section{The Generalized Preferential Attachment Model}\label{s:model}

\subsection{Model Description}\label{s:model.desc}

Preferential attachment is the rule by which new edges connect to existing nodes with probability proportional to their current degree:
\begin{equation}\label{eq:PA}
    p(u) \propto \frac{d(u)}{\sum_{v \in V} d(v)},
\end{equation}
where $d(v)$ is the degree (hyperdegree in the hypergraph setting) of node $v$. We extend preferential attachment to hypergraphs by introducing two random variables at each timestep: $Y_t$, governing size of the new hyperedge, and $X_t \leq Y_t$, governing how many of the new hyperedge's members are new nodes. Formally, at each timestep $t$ the model executes:
\begin{enumerate}
    \item Sample $Y_t$ from a distribution over the positive integers.
    \item Sample $X_t$ from a distribution over $\{0, 1, \ldots, Y_t\}$.
    \item Select $Y_t - X_t$ pre-existing nodes without replacement from the hypergraph according to rule \eqref{eq:PA}.
    \item Form a new hyperedge of size $Y_t$ containing the $X_t$ newly added nodes and the $Y_t - X_t$ preferentially selected pre-existing nodes.
\end{enumerate}
When $X_t = Y_t - 1$ deterministically, the model recovers the acyclic model of \cite{EvolvingHGModel}. Setting $X_t \in \{0, 1\}$ as a Bernoulli draw recovers the model of Avin et al.~\cite{PA_HG}. Our formulation generalizes both, allowing $X_t$ to take any value in $\{0, \ldots, Y_t\}$ and drawing both random variables from arbitrary distributions.

\subsection{Hyperdegree Distribution}\label{s:model.dist}

Let: $$p =\frac{ \lim_{t\to\infty} E[X_t] }{ \lim_{t\to\infty} E[Y_t]} = \frac{\overline{X}}{\overline{Y}}$$ the limiting ratio of expected new nodes to expected hyperedge size. We assume $\overline{X}$ and $\overline{Y}$ are both finite and strictly positive, and that $p \in (0,1)$ (see Remark~\ref{rem:boundary}).

\begin{theorem}\label{thm:main}
In the mean-field thermodynamic limit ($t \to \infty$), the stationary hyperdegree distribution of the generalized preferential attachment model is
\begin{equation}\label{eq:dist}
    P(d_H) = \frac{\Gamma(d_H)\,\Gamma\!\left(\tfrac{1}{1-p}+1\right)}
                  {(1-p)\,\Gamma\!\left(\tfrac{1}{1-p}+d_H+1\right)},
    \quad d_H \geq 1,
\end{equation}
and satisfies $P(d_H) \propto d_H^{-\left(\frac{1}{1-p}+1\right)}$ for $d_H \gg \max\!\left(1,\, \tfrac{1}{1-p}\right)$. The power-law exponent $\gamma = \tfrac{1}{1-p}+1 \in (2, \infty)$ depends only on $p$ and not on the shapes of the distributions of $X_t$ or $Y_t$ beyond their means.
\end{theorem}

The full derivation is given in Appendix~\ref{a:derivation}. The key steps are: (i) write an AME for the probability that a given node has hyperdegree $d_H$ at time $t$, replacing random quantities by their conditional expectations under the mean-field approximation; (ii) average over all nodes to obtain the time-dependent distribution $P(d_H, t)$; (iii) impose stationarity, appealing to the strong law of large numbers to justify the deterministic limits $D_t/t \to \overline{Y}$ and $|V_t|/t \to \overline{X}$ almost surely; (iv) solve the resulting recursion by induction on $d_H$ and express the result in terms of gamma functions.

The exponent $\gamma > 2$ guarantees a finite mean hyperdegree, a necessary condition for the preferential attachment mechanism to produce a well-defined large-$N$ limit.

\begin{corollary}\label{cor:ccdf}
The complementary cumulative distribution function satisfies
\begin{equation}\label{eq:ccdf}
    \mathrm{CCDF}(d_H) = \frac{\Gamma\!\left(\tfrac{1}{1-p}+1\right)
                               \,\Gamma(d_H+1)}
                              {\Gamma\!\left(\tfrac{1}{1-p}+d_H+1\right)},
\end{equation}
with asymptotic behavior $\mathrm{CCDF}(d_H) \propto d_H^{-1/(1-p)}$ for large $d_H$.
\end{corollary}

The proof is by induction on $d_H$; see Appendix~\ref{a:ccdf}.

\begin{remark}\label{rem:boundary}
The formula~\eqref{eq:dist} is derived for $p \in (0,1)$. The boundary case $p = 1$ is degenerate: when $X_t = Y_t$ at every step no pre-existing nodes are selected, so no preferential attachment occurs, each node appears in exactly one hyperedge with degree 1, and the stationary distribution is trivially $P(1) = 1$.
\end{remark}

\subsection{Relationship to Existing Models}\label{s:model.xi}

Setting $\xi = 1/(1-p)$, the distribution \eqref{eq:dist} can be written as
\begin{equation}
    P(d_H) = \frac{\xi}{\xi+1}\prod_{d=1}^{d_H-1}\frac{d}{\xi+1+d},
\end{equation}
matching the form derived by Avin et al.~\cite{PA_HG} for a model in which each hyperedge contains exactly one pre-existing node. Our result shows that this form is universal: the parameter $\xi = \overline{Y}/(\overline{Y}-\overline{X})$ (the ratio of average total hyperdegree contribution per step to the average number of pre-existing nodes being wired to per step) fully determines the degree distribution regardless of how $X_t$ and $Y_t$ are distributed.

\section{Nonlinear Preferential Attachment}\label{s:nlpa}

To investigate how the \emph{strength} of preferential attachment shapes simpliciality, we replace rule~\eqref{eq:PA} with the nonlinear variant
\cite{Gelation, NLPA}:
\begin{equation}\label{eq:NPA}
    p(u) \propto \frac{d(u)^\alpha}{\sum_{v\in V} d(v)^\alpha},
\end{equation}
where $\alpha \geq 0$ is the \emph{preferential exponent}. The limit $\alpha \to 0^+$ gives uniform random attachment; $\alpha = 1$ gives the standard linear rule; $\alpha \in (0,1)$ is sublinear; and $\alpha > 1$ is superlinear.

For $\alpha > 1$ the model undergoes a \emph{gelation transition} \cite{Gelation}: the attachment probability concentrates increasingly on the highest-degree node, and in the limit $N \to \infty$ a single hub accumulates a nonzero fraction of all hyperedges. In this regime the hub's neighborhood becomes an increasingly dense clique, which --- as we show in Section~\ref{s:results} --- disrupts downward closure and reduces simpliciality. The gelation transition at $\alpha = 1$ therefore defines a natural boundary: preferential attachment enforces simpliciality for $\alpha \leq 1$ and undermines it for $\alpha > 1$.

\section{Datasets and Fitting Methodology}\label{s:data}

\subsection{Datasets}\label{s:data.datasets}

We analyze timestamped hypergraph datasets drawn from the XGI-DATA repository \cite{XGI_data}, which compiles and standardizes datasets from the Austin Benson repository \cite{SimplicialClosure_LinkPrediction}, the SocioPatterns project, and additional sources. The datasets span four categories: collaboration networks (academic co-authorship and legislative co-sponsorship), communication networks (email threads), proximity contact networks (face-to-face contacts at 20-second resolution), and non-social tag networks. Summary statistics and empirical simpliciality measures are reported in Table~\ref{tab:network_stats}.

\subsection{Extracting $X_t$ and $Y_t$ from Data}\label{s:data.fitting}

A key feature of the generalized model is that its parameters can be estimated directly from any timestamped hypergraph without imposing distributional assumptions. Given a hypergraph in which each hyperedge carries a timestamp, we reconstruct the sequence $\{(X_t, Y_t)\}_{t=1}^T$ via a backward-stepping procedure: starting from the fully aggregated hypergraph $H_T$, remove the most recently added hyperedge $e_T$; record $Y_T = |e_T|$ and $X_T = $ the number of nodes left isolated by its removal (i.e., nodes whose only membership was $e_T$, and therefore first appeared when $e_T$ was added); repeat on the resulting $H_{T-1}$. The fitted ratio is then $\hat{p} = \overline{X}/\overline{Y}$ where $\overline{X}$ and $\overline{Y}$ are the sample means of the extracted sequences.

\begin{table}[htbp]
\centering
\small
\label{tab:network_stats}
\begin{tabular}{ll|rr|cc}
\toprule
\multicolumn{2}{c|}{\textbf{Dataset}}
  & \textbf{$|V|$} & \textbf{$|E|$}
  & \textbf{$\sigma_{SF}$} & \textbf{$\hat{p}$} \\
\midrule
\multirow{1}{*}{\textbf{Collaboration}}
 & Congress-Bills      &   1,718 & 127,127  & 0.019
   & 0.00078 \\
\midrule
\multirow{2}{*}{\textbf{Communication}}
 & Email-Enron         &     143 &  10,454  & 0.259 & 0.0055 \\
 & Email-EU            &     986 & 209,508  & 0.440 & 0.0020 \\
\midrule
\multirow{5}{*}{\textbf{Proximity}}
 & Contact-high-school    &   327 & 172,035  & 0.911 & 0.00098 \\
 & Contact-primary-school &   242 & 106,879  & 0.911 & 0.00095 \\
 & Hospital-Lyon          &    75 &  27,843  & 0.965 & 0.0014 \\
 & Malawi-Village         &    84 &  99,910  & 0.996 & 0.00050 \\
 & Science-Gallery        &   410 &  14,274  & 0.810 & 0.015 \\
\bottomrule
\end{tabular}
\caption{Dataset statistics, empirical simplicial fractions $\sigma_{SF}$, and fitted ratios $\hat{p}$ for all datasets. Entries marked in red are to be completed. NLPA comparison figures are shown for the eight datasets for which they are currently available.}
\end{table}

\subsection{Null Models}\label{s:data.nullmodels}

To determine whether the simpliciality of model-generated hypergraphs exceeds what is expected from the empirical structural statistics alone, we compare against four randomization null models that progressively preserve more of the original structure.

\begin{enumerate}
    \item \textbf{Random shuffling}: the nodes in each hyperedge are replaced with uniformly randomly selected nodes, preserving only hyperedge sizes.
    
    \item \textbf{Proportional shuffling}: nodes are replaced by nodes selected with probability proportional to their empirical hyperdegree, approximately preserving the degree distribution while destroying higher-order structure \cite{Hypergraph_Dissimilarity_Measures}.
    \item \textbf{Hyperdegree-preserving shuffling}: pairs of same-sized hyperedges are selected and one node is swapped between them, exactly preserving the hyperdegree sequence
    \cite{Hypergraph_Dissimilarity_Measures}.
    \item \textbf{Layer-preserving shuffling}: node labels are permuted within each layer (the set of hyperedges of a given size), preserving within-layer structure while randomizing cross-layer connections
    \cite{Encapsulation}.
\end{enumerate}

For each dataset, synthetic hypergraphs are generated under each null model and under the NLPA model across $\alpha \in [0, 2]$, computing $\sigma_{SF}$ for each and comparing against the empirical value.

\section{Results}\label{s:results}

\subsection{All Datasets Reside in the Small-$p$ Regime}\label{s:results.p}

Applying the backward-stepping procedure to all datasets, we find that $\hat{p}$ is extremely small across all categories, ranging from $\hat{p} \approx 0.0005$ for the most saturated proximity contact networks to $\hat{p} \approx 0.015$. Every dataset falls in the small-$p$ regime, arising naturally because the node sets of these systems are largely fixed: most hyperedges connect pre-existing participants ($X_t \approx 0$). This concentration of $X_t$ near zero is most extreme in proximity contact datasets, where hyperedges form among a bounded population, and in legislative datasets, where the same legislators co-sponsor many bills.

The power-law exponent predicted by Theorem~\ref{thm:main} satisfies $\gamma = 1/(1-\hat{p}) + 1 \approx 2$ in all cases, and the thermodynamic limit is approached slowly due to finite-size effects in the heavy tail (see Appendix~\ref{a:validation}). We therefore do not attempt to validate Theorem~\ref{thm:main} directly against empirical degree sequences in these datasets; the validation is instead carried out against controlled simulations in Appendix~\ref{a:validation}, where both $\hat{p}$ and $N$ can be controlled independently.

Figure~\ref{fig:sf_vs_p} shows the simplicial fraction of model-generated hypergraphs as a function of $p$, for the generalized model under linear preferential attachment ($\alpha = 1$). The simplicial fraction decreases rapidly and reaches effectively zero before $p = 0.25$ across all values of $E[Y_t]$ considered, confirming that non-trivial simplicial structure in the model is confined to the small-$p$ regime, which encompasses all empirical datasets studied.

\subsection{Preferential Attachment as a Simpliciality-Enforcing Mechanism}\label{s:results.nlpa}

\begin{figure}[htbp]
\includegraphics[width=\textwidth]{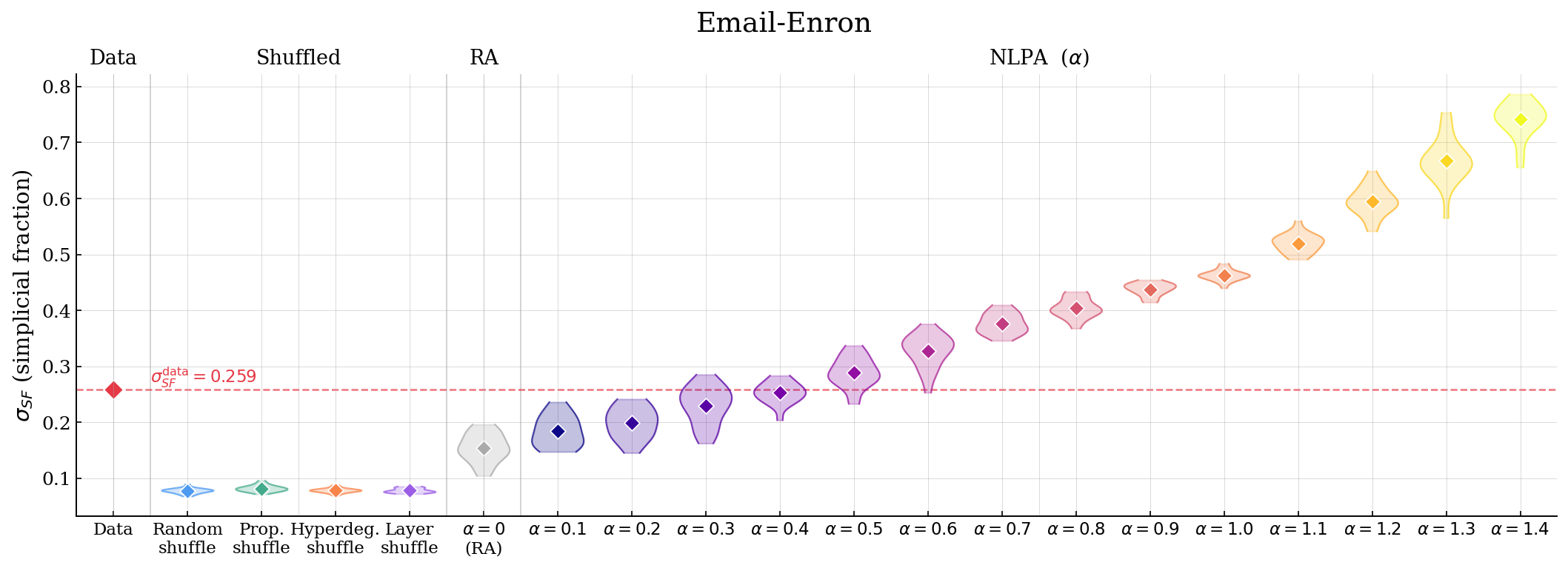}
\centering
\caption{Simplicial fraction $\sigma_{SF}$ for Email-Enron ($\sigma_{SF}^{\mathrm{data}} = 0.259$, dashed red line) as a function of the preferential exponent $\alpha$, compared against the four null models and the NLPA model. All null models produce $\sigma_{SF} \approx 0$, confirming that the observed simpliciality is not explained by the hyperedge size or degree distributions alone. The NLPA model matches the empirical value at approximately $\alpha \approx 0.4$.}
\label{fig:nlpa_enron}
\end{figure}

\begin{figure}[htbp]
\includegraphics[width=\textwidth]{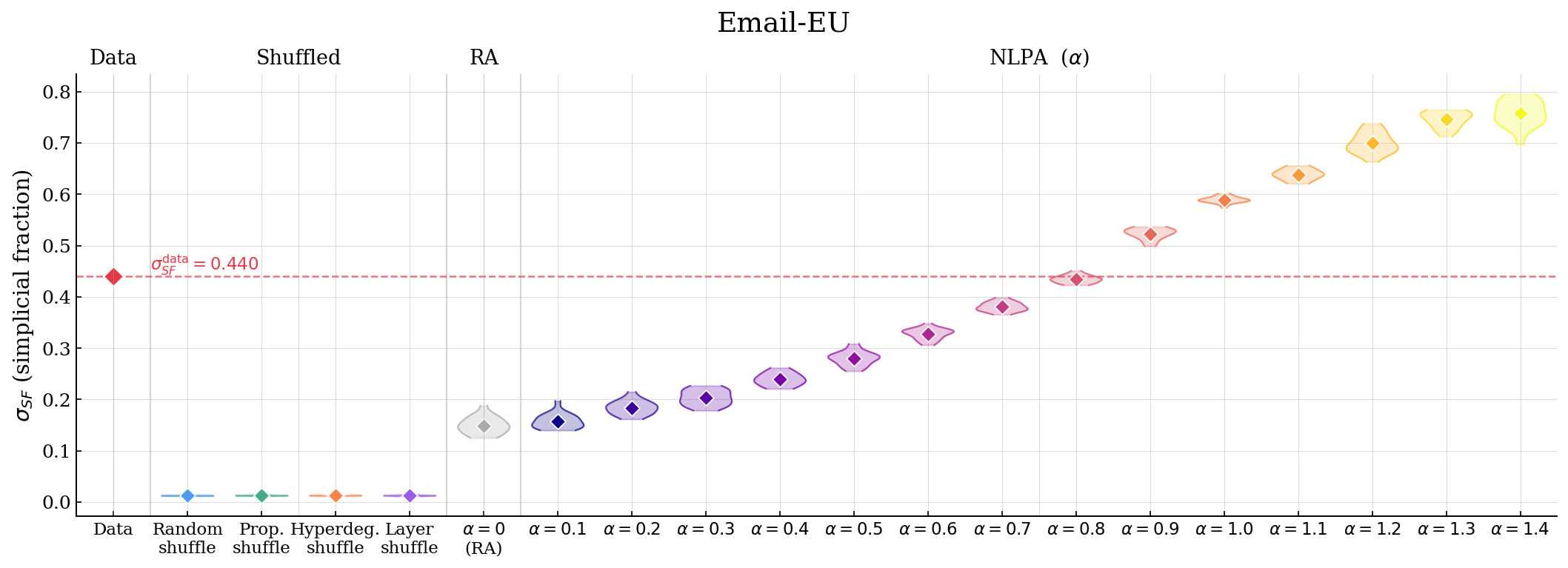}
\centering
\caption{Simplicial fraction $\sigma_{SF}$ for Email-EU ($\sigma_{SF}^{\mathrm{data}} = 0.440$, dashed red line). As in Figure~\ref{fig:nlpa_enron}, null models produce near-zero $\sigma_{SF}$, confirming attachment-driven simpliciality. The NLPA model matches the empirical value at approximately $\alpha \approx 0.8$.}
\label{fig:nlpa_eu}
\end{figure}

\begin{figure}[htbp]
\includegraphics[width=\textwidth]{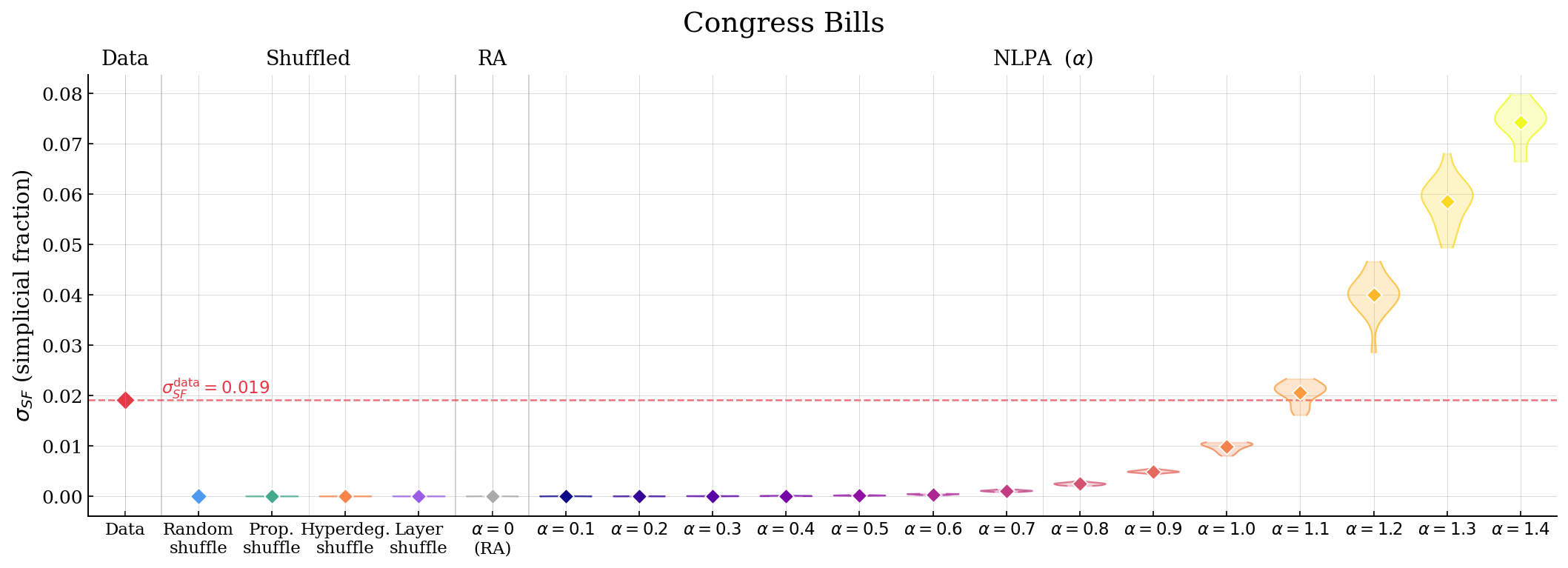}
\centering
\caption{Simplicial fraction $\sigma_{SF}$ for Congress-Bills ($\sigma_{SF}^{\mathrm{data}} = 0.019$, dashed red line). Despite the low absolute simpliciality, all null models and random attachment produce $\sigma_{SF} \approx 0$, and the empirical value is matched only at $\alpha \approx 1.1$, in the superlinear regime just beyond the gelation threshold.}
\label{fig:nlpa_congress}
\end{figure}

\begin{figure}[htbp]
\includegraphics[width=\textwidth]{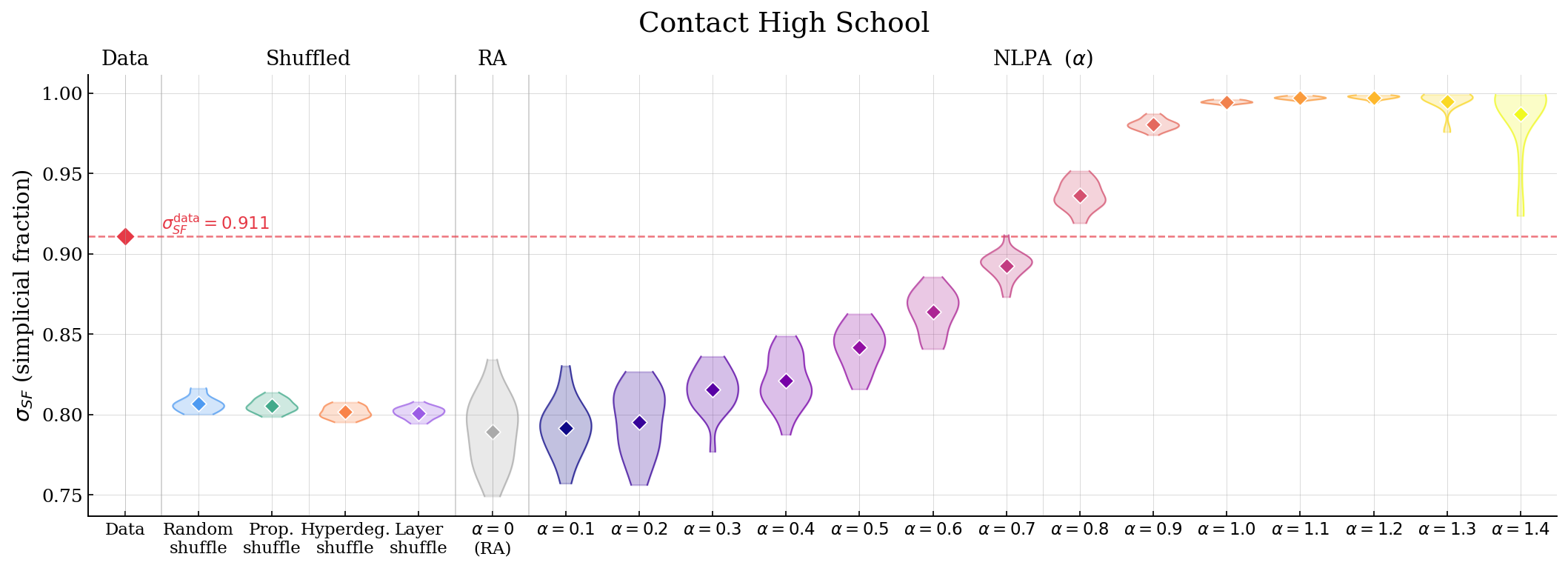}
\centering
\caption{Simplicial fraction $\sigma_{SF}$ for Contact-High-School ($\sigma_{SF}^{\mathrm{data}} = 0.911$, dashed red line). Null models produce $\sigma_{SF} \approx 0.80$--$0.81$, indicating that both the hyperedge size distribution and the attachment mechanism contribute to simpliciality. The NLPA model matches the empirical value at $\alpha$ between $0.7$ and $0.8$.}
\label{fig:nlpa_hs}
\end{figure}

\begin{figure}[htbp]
\includegraphics[width=\textwidth]{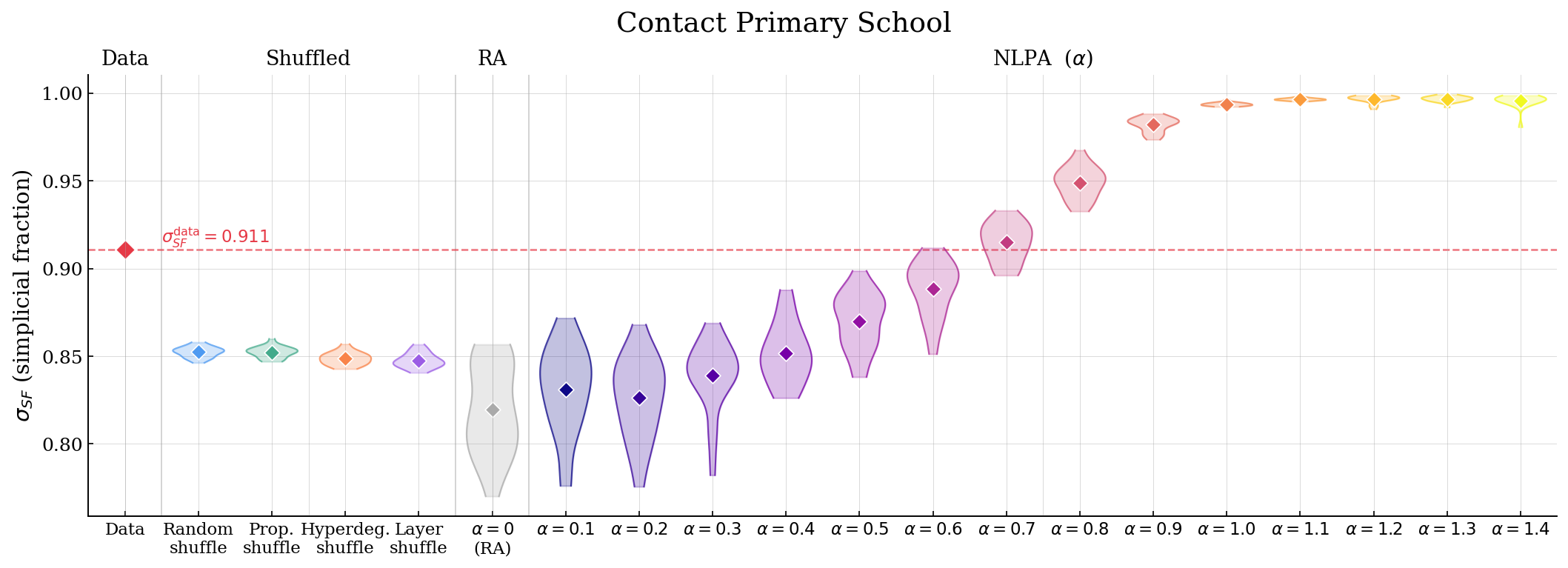}
\centering
\caption{Simplicial fraction $\sigma_{SF}$ for Contact-Primary-School ($\sigma_{SF}^{\mathrm{data}} = 0.911$, dashed red line). Null models produce $\sigma_{SF} \approx 0.85$ and the NLPA model matches the empirical value at approximately $\alpha \approx 0.7$.}
\label{fig:nlpa_ps}
\end{figure}

\begin{figure}[htbp]
\includegraphics[width=\textwidth]{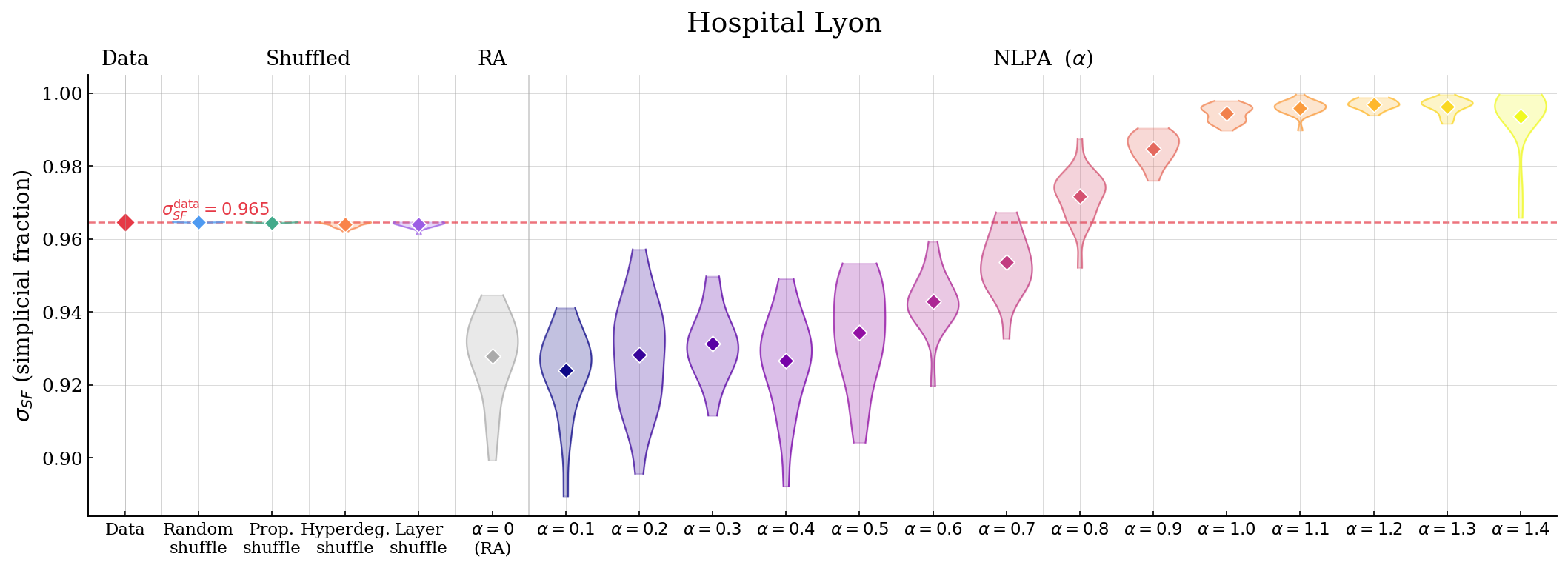}
\centering
\caption{Simplicial fraction $\sigma_{SF}$ for Hospital-Lyon ($\sigma_{SF}^{\mathrm{data}} = 0.965$, dashed red line). All four null models produce $\sigma_{SF} \approx 0.965$, matching the empirical value, while random attachment ($\alpha = 0$) yields only $\approx 0.929$. The near-maximal simpliciality is primarily encoded in the hyperedge size and degree distributions, with preferential attachment providing a secondary reinforcing role. The NLPA model matches the empirical value at $\alpha \approx 0.75$.}
\label{fig:nlpa_hospital}
\end{figure}

\begin{figure}[htbp]
\includegraphics[width=\textwidth]{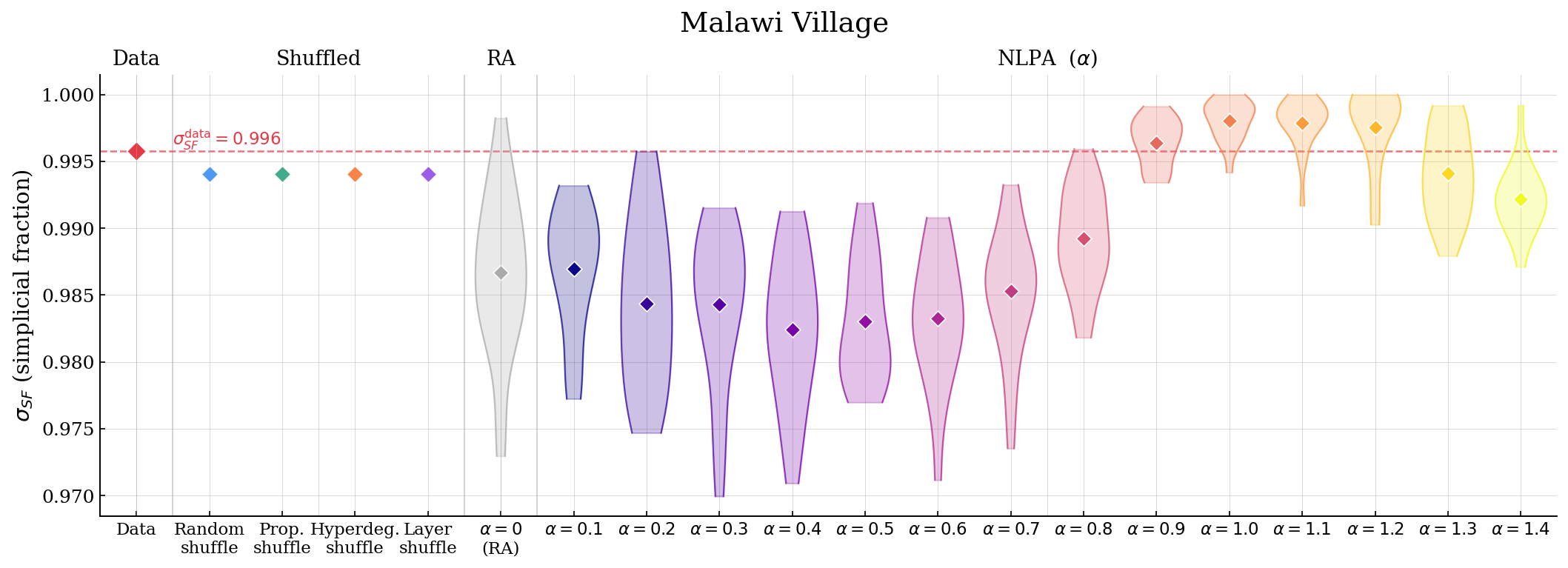}
\centering
\caption{Simplicial fraction $\sigma_{SF}$ for Malawi-Village ($\sigma_{SF}^{\mathrm{data}} = 0.996$, dashed red line). As with Hospital-Lyon, all null models reproduce the empirical simpliciality ($\sigma_{SF} \approx 0.994$), confirming distribution-driven simpliciality. The NLPA model matches the empirical value at $\alpha \approx 0.9$.}
\label{fig:nlpa_malawi}
\end{figure}

\begin{figure}[htbp]
\includegraphics[width=\textwidth]{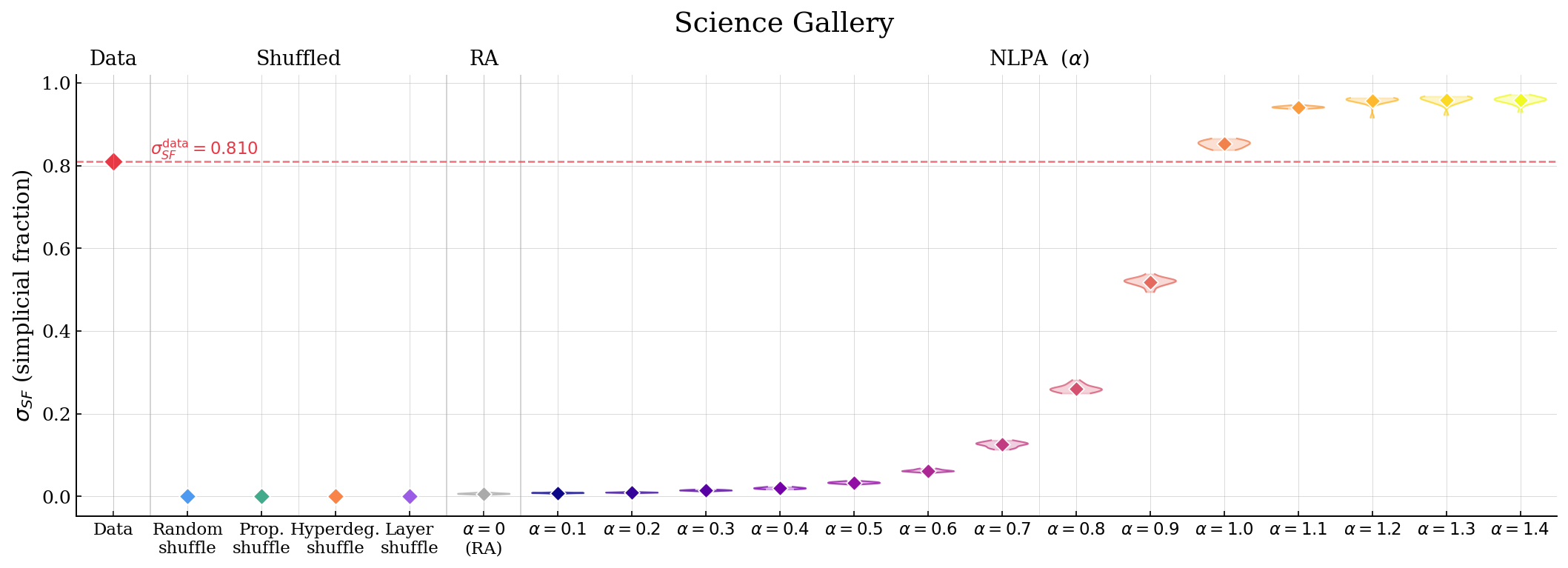}
\centering
\caption{Simplicial fraction $\sigma_{SF}$ for Science-Gallery ($\sigma_{SF}^{\mathrm{data}} = 0.810$, dashed red line). Null models and random attachment both produce $\sigma_{SF} \approx 0$, and the NLPA model rises steeply near $\alpha \approx 1.0$ to match the empirical value, consistent with attachment-driven simpliciality at approximately the linear preferential attachment threshold.}
\label{fig:nlpa_science}
\end{figure}

Figures~\ref{fig:nlpa_enron}--\ref{fig:nlpa_science} show $\sigma_{SF}$ of
NLPA-generated hypergraphs as a function of $\alpha \in [0, 2]$, compared against the four null models and the empirical value (dashed red line) for each of the eight datasets with available results. A consistent finding across all datasets is that $\sigma_{SF}$ increases monotonically with $\alpha$ for $\alpha \leq 1$, and decreases for $\alpha > 1$. This establishes preferential attachment as a simpliciality-enforcing mechanism in the sublinear and linear regimes, with the gelation transition at $\alpha = 1$ marking the boundary beyond which stronger attachment begins to undermine simpliciality.

The mechanism is straightforward: below the gelation threshold, increasing $\alpha$ causes pre-existing high-degree nodes to be selected more consistently, concentrating edges around nodes whose neighborhoods become dense enough that many sub-hyperedges are shared, satisfying downward closure. Beyond gelation, the dominant hub accumulates hyperedges so aggressively that newly added hyperedges connect primarily to the hub and its immediate neighbors, creating large hyperedges whose sub-hyperedges are absent from the hypergraph.

The datasets divide into two categories based on the relationship between null model simpliciality and the empirical value.

\textbf{Attachment-driven simpliciality.} For Email-Enron (Figure~\ref{fig:nlpa_enron}), Email-EU (Figure~\ref{fig:nlpa_eu}), and Congress-Bills (Figure~\ref{fig:nlpa_congress}), all four null models produce $\sigma_{SF} \approx 0$, substantially below the empirical values of $0.259$, $0.440$, and $0.019$ respectively. The hyperedge size and degree distributions alone are insufficient to account for the observed simpliciality; the attachment mechanism is the primary driver. The empirical $\sigma_{SF}$ is matched at approximately $\alpha \approx 0.4$ (Email-Enron), $\alpha \approx 0.8$ (Email-EU), and $\alpha \approx 1.1$ (Congress-Bills). The last of these lies in the superlinear regime, indicating that the modest but nonzero simpliciality of legislative co-sponsorship networks is consistent with a stronger-than-linear attachment mechanism. Contact-High-School (Figure~\ref{fig:nlpa_hs}), Contact-Primary-School (Figure~\ref{fig:nlpa_ps}), and Science-Gallery (Figure~\ref{fig:nlpa_science}) occupy an intermediate position: null models produce elevated but sub-empirical $\sigma_{SF}$, and the attachment mechanism contributes an additional increment matched at $\alpha$ between $0.7$ and $1.0$.

\textbf{Distribution-driven simpliciality.} For Hospital-Lyon (Figure~\ref{fig:nlpa_hospital}) and Malawi-Village (Figure~\ref{fig:nlpa_malawi}), all four null models produce $\sigma_{SF}$ values comparable to the empirical data ($\approx 0.965$ and $\approx 0.994$ respectively), while random attachment ($\alpha = 0$) falls below. In these saturated networks, the near-maximal simpliciality is largely determined by the hyperedge size and degree distributions themselves. Preferential attachment plays a secondary but consistent reinforcing role, with the empirical value matched at $\alpha \approx 0.75$ (Hospital-Lyon) and $\alpha \approx 0.9$ (Malawi-Village).

Across all datasets and categories, the monotone dependence of $\sigma_{SF}$ on $\alpha$ in the sublinear regime confirms that preferential attachment is a structurally coherent simpliciality-enforcing mechanism, regardless of domain or scale.

\section{Discussion}\label{s:discussion}

The central analytical contribution of this paper is the universality result of Theorem~\ref{thm:main}: the power-law exponent of the stationary hyperdegree distribution depends on $X_t$ and $Y_t$ only through their ratio $p = \overline{X}/\overline{Y}$. Any two choices of distributions sharing the same ratio produce identical degree distributions, regardless of how the individual sizes are distributed.

\begin{figure}[htbp]
\includegraphics[width=.6\textwidth]{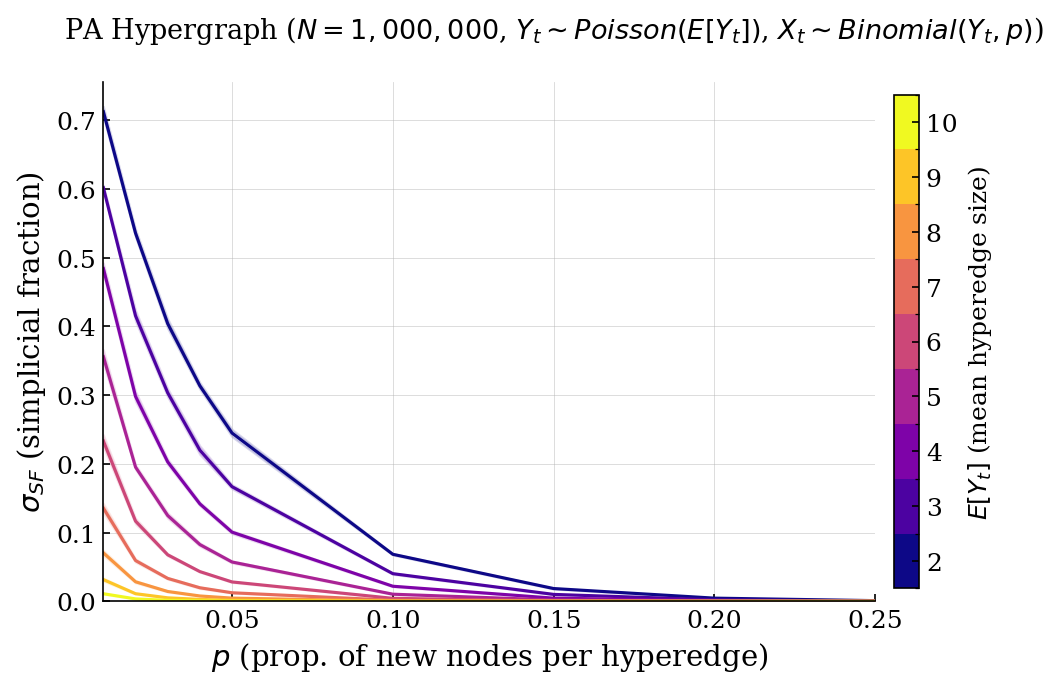}
\centering
\caption{Simplicial fraction $\sigma_{SF}$ of generalized preferential attachment hypergraphs ($N = 1{,}000{,}000$) as a function of $p$, for $E[Y_t] \in \{2, 3, \ldots, 10\}$ under linear preferential attachment ($\alpha = 1$). For all values of $E[Y_t]$, $\sigma_{SF}$ drops to effectively zero before $p = 0.25$, confirming that non-trivial simplicial structure is confined to the small-$p$ regime. All empirical datasets studied have $\hat{p} \ll 0.25$.}
\label{fig:sf_vs_p}
\end{figure}

The observation that all studied datasets reside in the small-$p$ regime has two important consequences. First, the analytical degree distribution converges slowly to the thermodynamic limit, so the power-law character is most visible at large hyperdegree and large system size. Second, the simplicial fraction of model-generated hypergraphs is non-trivially sensitive to the attachment mechanism in this regime, as shown in Figure~\ref{fig:sf_vs_p}: small $p$ is precisely where preferential attachment has the greatest potential to enforce simpliciality.

The NLPA results reveal that different types of higher-order interaction systems have qualitatively different drivers of simpliciality. In sparse, open-membership systems such as email and legislative networks, preferential attachment is the primary mechanism and the fitted $\alpha$ provides a natural summary statistic of simplicial excess beyond what structural statistics predict. In dense, closed-population proximity networks, the degree and size distributions account for most of the simpliciality and the attachment mechanism plays a secondary role. This suggests that simpliciality has different origins in different classes of higher-order systems, and that care is needed before applying topological tools to hypergraph data without first assessing the degree of simpliciality and its origin.

\section{Conclusion \& Future Work}\label{s:conclusion}

We have introduced a generalized preferential attachment hypergraph model and established three results. First, using a mean-field AME approach, the stationary hyperdegree distribution follows a power law whose exponent depends only on $p = E[X_t]/E[Y_t]$, providing an analytically tractable framework for understanding degree heterogeneity in higher-order systems. Second, the model parameters $X_t$ and $Y_t$ can be estimated directly from any timestamped hypergraph dataset, grounding the model in empirical network statistics. Third, across diverse real-world datasets, the simplicial fraction increases monotonically with the preferential attachment exponent up to the gelation transition, establishing preferential attachment as a simpliciality-enforcing mechanism.

Several directions for future work follow naturally. The small-$p$ regime warrants dedicated finite-size analysis, as convergence to the thermodynamic limit is slow and the power law is not visible at typical empirical system sizes. The backward-stepping fitting procedure could be extended to handle simultaneous hyperedge formation and non-stationary node arrival rates. Whether simpliciality enforcement persists under fitness-based \cite{BA_Fitness} or local-world \cite{PA_HG_LocalWorld} attachment kernels remains an open question. Finally, comparing preferential attachment against other generative mechanisms such as edge-copying \cite{NoisyCopyModel} on the same empirical datasets would help determine the extent to which preferential attachment uniquely accounts for observed simpliciality.


\subsubsection*{Data Availability}
All datasets are available in the XGI-DATA repository \cite{XGI_data}. Code implementing the generalized and nonlinear preferential attachment models, simpliciality computations, and figure generation is openly available at \texttt{https://github.com/JasonLaRuez/JasonLaRuez\_Dissertation}.

\printbibliography

\appendix

\section{Derivation of the Hyperdegree Distribution}\label{a:derivation}

\noindent We derive the hyperdegree distribution for the generalized preferential attachment model in which hyperedge sizes $Y_t$ and new-node counts $X_t \leq Y_t$ are drawn from arbitrary distributions at each timestep.

\paragraph{Setup and notation.} Let $d_H$ be a particular value of the hyperdegree. Let $v_{(i,j)}$ denote the $j$-th node added at time $i$ ($j \leq X_i$), and let $P(d_H, v_{(i,j)}, t)$ be the probability that $v_{(i,j)}$ has hyperdegree $d_H$ at time $t$. Let $|V_t| = \sum_{i=1}^t X_i$ and $D_t = \sum_{v \in V_t} d(v) = \sum_{i=1}^t Y_i$ be the total node count and total hyperdegree at time $t$, respectively. Note that both $|V_t|$ and $D_t$ are sums of i.i.d.\ random variables; by the strong law of large numbers (SLLN),
\begin{equation}\label{eq:slln}
    \frac{|V_t|}{t} \to \overline{X}, \quad
    \frac{D_t}{t} \to \overline{Y} \quad \text{almost surely as } t \to \infty,
\end{equation}
so that $|V_t|/D_t \to \overline{X}/\overline{Y} = p$ almost surely.

\paragraph{The approximate master equation.} At time $t$, the hyperedge added consists of $X_t$ new nodes and $Y_t - X_t$ pre-existing nodes selected preferentially. We make two standard mean-field approximations. First, we replace the random quantity $Y_t - X_t$ by its conditional expectation $E[Y_t - X_t]$. Second, we treat the $Y_t - X_t$ preferential selections as independent Bernoulli events, each succeeding for node $v_{(i,j)}$ with probability $d_H/D_t$. Under these approximations:
\begin{equation}\label{eq:ame}
\begin{split}
    P(d_H, v_{(i,j)}, t+1) &=
    \frac{(d_H-1)(Y_t-X_t)}{D_t}\,P(d_H-1, v_{(i,j)}, t) \\
    &\quad + \left(1 - \frac{d_H(Y_t-X_t)}{D_t}\right)
    P(d_H, v_{(i,j)}, t).
\end{split}
\end{equation}

\paragraph{Averaging over nodes.} The time-dependent distribution is
\begin{equation}\label{eq:timedep}
    P(d_H, t) = \frac{1}{|V_t|}\sum_{i=1}^{t}\sum_{j=1}^{X_i}
    P(d_H, v_{(i,j)}, t).
\end{equation}
Multiplying both sides of \eqref{eq:timedep} at time $t+1$ by $|V_{t+1}|$ and separating the $i = t+1$ term (nodes added at time $t+1$ have hyperdegree 1 by construction, contributing $X_{t+1} \delta_1(d_H)$):
\begin{align*}
    |V_{t+1}| P(d_H,t+1)
    &= (d_H-1)(Y_t-X_t)\frac{|V_t|}{D_t}\,P(d_H-1,t) \\
    &\quad + \left(1 - \frac{d_H(Y_t-X_t)}{D_t}\right)|V_t|\,P(d_H,t)
    + X_{t+1}\delta_1(d_H).
\end{align*}
where $\delta_i(j)$ denotes the Kronecker delta  function. Rearranging and using $|V_{t+1}| = |V_t| + X_{t+1}$:
\begin{align*}
    X_{t+1}\,P(d_H,t+1) + |V_t|\,\frac{P(d_H,t+1)-P(d_H,t)}{1}
    &= (d_H-1)(Y_t-X_t)\frac{|V_t|}{D_t}\,P(d_H-1,t) \\
    &\quad - d_H(Y_t-X_t)\frac{|V_t|}{D_t}\,P(d_H,t)
    + X_{t+1}\delta_1(d_H).
\end{align*}
We apply the continuous-time approximation $(P(d_H, t+1) - P(d_H, t))/((t+1)-t) \approx \partial P/\partial t$, valid near stationarity where $P(d_H, t)$ varies on timescales $\gg 1$.

\paragraph{Stationary limit.} Taking $t \to \infty$, stationarity requires $\partial P/\partial t = o(1/t)$, so the $t\,\partial P/\partial t$ term vanishes. Applying the SLLN limits \eqref{eq:slln} and replacing random quantities by their limiting expectations:
\begin{align*}
    X_{t+1} &\to \overline{X} = p\overline{Y}, \\
    Y_t - X_t &\to \overline{Y} - \overline{X} = (1-p)\overline{Y}, \\
    |V_t|/D_t &\to p.
\end{align*}
The stationary equation becomes:
\begin{equation}
    p\overline{Y}\,P(d_H) = (d_H-1)(1-p)\overline{Y}p\,P(d_H-1)
    - d_H(1-p)\overline{Y}p\,P(d_H)
    + p\overline{Y}\delta_1(d_H).
\end{equation}
Dividing by $p\overline{Y}$ and rearranging:
\begin{equation}\label{eq:recursion}
    P(d_H) = \frac{(d_H-1)(1-p)\,P(d_H-1) + \mathbf{1}[d_H=1]}
                  {1 + d_H(1-p)}.
\end{equation}

\paragraph{Solving the recursion.} We claim
\begin{equation}\label{eq:inductform}
    P(d_H) = \frac{(d_H-1)!\,(1-p)^{d_H-1}}
                  {\displaystyle\prod_{d=1}^{d_H}\bigl(1+d(1-p)\bigr)},
    \quad d_H \geq 1.
\end{equation}
\textit{Base case} ($d_H = 1$): recursion \eqref{eq:recursion} gives $P(1) = 1/[1+(1-p)] = 1/(2-p)$, which equals \eqref{eq:inductform}. $\checkmark$ \\

\noindent \textit{Inductive step}: Assume \eqref{eq:inductform} holds for $d_H - 1$.
Then:
\begin{align*}
    P(d_H)
    &= \frac{(d_H-1)(1-p)}{1+d_H(1-p)} \cdot
       \frac{(d_H-2)!\,(1-p)^{d_H-2}}
            {\prod_{d=1}^{d_H-1}(1+d(1-p))} \\
    &= \frac{(d_H-1)!\,(1-p)^{d_H-1}}
            {\prod_{d=1}^{d_H}(1+d(1-p))}, \qquad \checkmark
\end{align*}
completing the induction.

\paragraph{Gamma function form.} Using $\Gamma(z+1) = z\Gamma(z)$ repeatedly:
\begin{equation*}
    \Gamma\!\left(\tfrac{1}{1-p}+d_H+1\right)
    = \Gamma\!\left(\tfrac{1}{1-p}+1\right)
      \prod_{d=1}^{d_H}\!\left(\tfrac{1}{1-p}+d\right)
    = \frac{\Gamma\!\left(\tfrac{1}{1-p}+1\right)}{(1-p)^{d_H}}
      \prod_{d=1}^{d_H}\!\bigl(1+d(1-p)\bigr).
\end{equation*}
Substituting into \eqref{eq:inductform} and using $(d_H-1)! = \Gamma(d_H)$:
\begin{equation}\label{eq:gammadist}
    \boxed{P(d_H) = \frac{\Gamma(d_H)\,\Gamma\!\left(\tfrac{1}{1-p}+1\right)}
                         {(1-p)\,\Gamma\!\left(\tfrac{1}{1-p}+d_H+1\right)}},
\end{equation}
which is the formula of Theorem~\ref{thm:main}.


\paragraph{Power-law tail.} Applying Stirling's approximation $\Gamma(z) \approx \sqrt{2\pi}\,z^{z+1/2}e^{-z}$ for large $z$, valid here for $d_H \gg \max(1,\, 1/(1-p))$:
\begin{align*}
    P(d_H)
    &\propto \frac{\Gamma(d_H)}{\Gamma\!\left(\tfrac{1}{1-p}+d_H+1\right)}
    \approx \frac{d_H^{d_H+1/2}\,e^{-d_H}}
               {\left(\tfrac{1}{1-p}+d_H+1\right)^{\frac{1}{1-p}+d_H+3/2}
                e^{-\left(\frac{1}{1-p}+d_H+1\right)}}.
\end{align*}
Using $(d_H/(\frac{1}{1-p}+d_H+1))^{d_H+1/2} \approx e^{-(\frac{1}{1-p}+1)}$ and $(\frac{1}{1-p}+d_H+1)^{\frac{1}{1-p}+1} \approx d_H^{\frac{1}{1-p}+1}$ for $d_H \gg 1/(1-p)$:
\begin{equation*}
    P(d_H) \propto d_H^{-\left(\frac{1}{1-p}+1\right)},
\end{equation*}
confirming the power law with exponent $\gamma = 1/(1-p)+1$.

\section{Proof of Corollary~\ref{cor:ccdf}: CCDF}\label{a:ccdf}

We claim $\mathrm{CCDF}(d_H) = \Gamma\!\left(\tfrac{1}{1-p}+1\right) \Gamma(d_H+1)/\Gamma\!\left(\tfrac{1}{1-p}+d_H+1\right)$.

\textit{Base case} ($d_H = 0$): $\mathrm{CCDF}(0) = 1 = \Gamma(\tfrac{1}{1-p}+1)\Gamma(1)/\Gamma(\tfrac{1}{1-p}+1)$. $\checkmark$

\textit{Inductive step}: Assume the formula holds for $d_H$. Then $\mathrm{CCDF}(d_H+1) = \mathrm{CCDF}(d_H) - P(d_H+1)$. Substituting:
\begin{align*}
    &= \frac{\Gamma\!\left(\tfrac{1}{1-p}+1\right)\Gamma(d_H+1)}
            {\Gamma\!\left(\tfrac{1}{1-p}+d_H+1\right)}
    - \frac{\Gamma(d_H+1)\,\Gamma\!\left(\tfrac{1}{1-p}+1\right)}
           {(1-p)\,\Gamma\!\left(\tfrac{1}{1-p}+d_H+2\right)} \\
    &= \Gamma\!\left(\tfrac{1}{1-p}+1\right)\Gamma(d_H+1)
       \cdot \frac{\tfrac{1}{1-p}+d_H+1 - \tfrac{1}{1-p}}
                  {\Gamma\!\left(\tfrac{1}{1-p}+d_H+2\right)} \\
    &= \frac{\Gamma\!\left(\tfrac{1}{1-p}+1\right)\,\Gamma(d_H+2)}
            {\Gamma\!\left(\tfrac{1}{1-p}+d_H+2\right)}, \qquad \checkmark
\end{align*}
where the final step uses $\Gamma(d_H+2) = (d_H+1)\Gamma(d_H+1)$.

\section{Simulation Validation}\label{a:validation}

\begin{figure}[htbp]
\includegraphics[width=\textwidth]{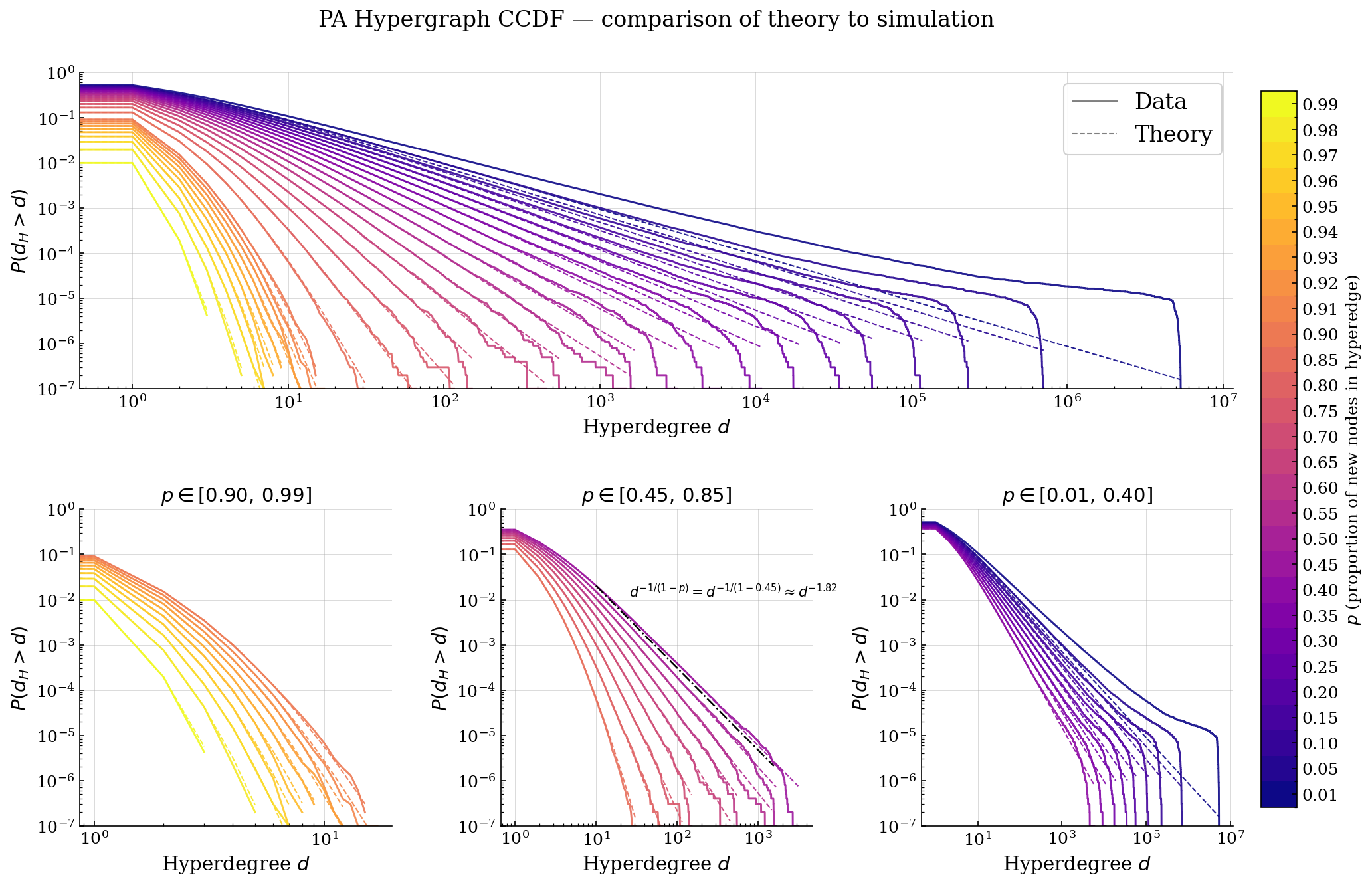}
\centering
\caption{Empirical CCDF of the hyperdegree distribution (solid lines) versus the theoretical prediction of Corollary~\ref{cor:ccdf} (dashed lines), for $Y_t \sim \mathrm{Poisson}(5)$ and $X_t \sim \mathrm{Binomial}(Y_t, p)$, across $p \in [0.01, 0.99]$. Top panel: all values of $p$. Bottom panels: three subranges for clarity. For large and intermediate $p$ the agreement is excellent across the full range of hyperdegrees. For small $p$, the empirical tail is heavier than the theoretical prediction: when $p$ is small, the same high-degree nodes are selected repeatedly and accumulate hyperdegree at a rate exceeding the mean-field prediction, a finite-size effect that becomes increasingly pronounced as $p \to 0$.}
\label{fig:ccdf}
\end{figure}

\noindent We validate Theorem~\ref{thm:main} against simulations of the generalized model with $Y_t \sim \mathrm{Poisson}(5)$ and $X_t \sim \mathrm{Binomial}(Y_t, p)$. The comparison between empirical and theoretical CCDFs across $p \in [0.01, 0.99]$ is shown in Figure~\ref{fig:ccdf}.

\end{document}